\documentclass[11pt,a4paper]{article}

\usepackage{amsmath,amssymb,amsthm,amscd}

\usepackage{array}
\newcolumntype{C}{>{\centering\arraybackslash}p{1.9ex}}

\usepackage{tikz}
\newcommand*\circled[1]{\tikz[baseline=(char.base)]{
            \node[shape=circle,draw,inner sep=.5pt] (char) {#1};}}

\newtheorem{lemma}{Lemma}

\title{On the b-chromatic number of the Cartesian product of two
complete graphs}

\author{Fr\'ed\'eric Maffray\thanks{CNRS, Laboratoire G-SCOP,
University of Grenoble, France.  Partially supported by ANR Project
Stint ANR-13-BS02-0007.  E-mail:
frederic.maffray@grenoble-inp.fr}%
\and%
Artur Mesquita Barbosa\thanks{Universidade Federal do Cear\'a,
Fortaleza, CE, Brazil.}}

\date{\today} 

\begin{document}
\maketitle

\section{Introduction}

A b-coloring of a graph $G$ is a proper coloring such that every color
class contains a vertex that has neighbors in all other classes
\cite{IM}.  Such a vertex will be called a \emph{b-vertex}.  The
b-chromatic number of $G$ is the largest integer $k$ such that $G$
admits a b-coloring with $k$ colors.  Let $\chi(G)$ be the usual
chromatic number of $G$, and let $\Delta(G)$ be the maximum degree in
$G$.  Let $t(G)$ be the largest integer such that $G$ has at least $t$
vertices of degree at least $t-1$.  If a graph $G$ has a b-coloring
with $k$ colors then $G$ has at least $k$ vertices of degree at least
$k-1$ (the b-vertices); hence the following chain of inequalities
holds for every graph:
$$\chi(G) \le b(G) \le t(G) \le \Delta(G)+1.$$

The Cartesian product of two graphs $G_1$ and $G_2$, denoted by
$G_1\square G_2$, is the graph with vertex-set $V(G_1) \times V (G_2)$
where two vertices $(u_1, u_2)$ and $(v_1, v_2)$ are adjacent if and
only if either $u_1 = v_1$ and $u_2v_2\in E(G_2)$ or $u_2=v_2$ and
$u_1v_1\in E(G_1)$.  Faik \cite{F} and Javadi and Omoomi \cite{JavOm}
studied b-colorings of the Cartesian product of two graphs for various
types of factor graphs.  Let $K_n$, $P_n$ and $C_n$ respectively
denote the complete graph, path and cycle on $n$ vertices.  In the
case of the Cartesian product of two copies of $K_n$, it is easy to
see that $\chi(K_n\square K_n)=n$ and $t(K_n\square K_n)=2n-1$.  Hence
$n\le b(K_n\square K_n)\le 2n-1$.  Javadi and Omoomi \cite{JavOm}
improved these bounds significantly as follows.
\begin{itemize}
\item
For all $n \ge 5$, we have $b(K_n\square K_n) \ge 2n - 3$.
\item
For all $n \ge 2$, we have $b(K_n\square K_n) \le 2n - 2$.
\end{itemize}
It follows that for all $n\ge 5$ we have $b(K_n\square K_n)
\in\{{2n-3}, {2n-2}\}$.  One would like to know which of these two
values is the right one.  Javadi and Omoomi \cite{JavOm} conjectured
that $b(K_n\square K_n) = 2n - 3$ for all $n\ge 5$.  We give
counterexamples to this conjecture: for $n=5$, so $b(K_5\square
K_5)=8$; for $n=6$, so $b(K_6\square K_6)=10$; and for $n=7$, so
$b(K_7\square K_7)=12$.  The counterexamples were found using a
computer search.  We explain below the ideas that helped speed up the
search; unfortunately these ideas were not sufficient to solve any
case with $n\ge 8$, because the number of possibilities to examine is
again too high.

\section{Results}

Let $G=K_5\square K_5$.  Let $G$ have vertex-set $\{v_{i,j} \mid i, j
= 1, \ldots, 5\}$ and edge-set $\{v_{i,j}v_{i,k}\mid i, j, k = 1,
\ldots, 5, j\neq k\} \cup \{v_{i,j}v_{h,j}\mid h, i, j = 1, \ldots, 5,
h\neq i\}$.  For each $i=1,\ldots,5$, the set $\{v_{i,1}, \ldots,
v_{i,5}\}$ is called \emph{row~$i$} of $G$, and the set $\{v_{1,i},
\ldots, v_{5,i}\}$ is called \emph{column~$i$} of $G$.  Assume that
$G$ admits a b-coloring $f$ with $8$ colors.  For each $i=1,\ldots,
5$, let $\beta_i$ be the number of b-vertices in row~$i$, and let
$\beta'_i$ be the number of b-vertices in column~$i$.

\begin{lemma}\label{lem:beta}
For each $i=1,\ldots,5$ we have $\beta_i\le 2$ and $\beta'_i\le 2$.
\end{lemma}
The proof of this lemma is given in the Appendix below.

\medskip

For $c=1, \ldots, 8$, pick one b-vertex $b_c$ of color $c$, and let
$B=\{b_1, \ldots, b_8\}$.  Let $H$ be the graph with vertex-set
$\{r_1, \ldots, r_5,$ $\ell_1, \ldots,$ $\ell_5\}$ such that $\{r_1,
\ldots,$ $r_5\}$ and $\{\ell_1, \ldots, \ell_5\}$ are stable sets and
$r_i\ell_j$ is an edge if and only if $v_{i,j}\in B$.  We call $H$ the
\emph{pattern} graph of $B$.  So $H$ is a bipartite graph with eight
edges, and by Lemma~\ref{lem:beta} every vertex of $H$ has degree at
most $2$.  Therefore every component of $H$ is a path or an even
cycle, and $H$ is one of the twelve graphs displayed in the following
table (where $+$ denotes the disjoint union and $2X$ denotes $X+X$):
\begin{center}
    \begin{tabular}{|l|l|l|l|}
	\hline
$C_8 + 2P_1$ & $C_6 + 2P_2$ & $C_6 + P_3 + P_1$ & $2C_4 + 2P_1$ \\ 
\hline
$C_4 + P_5 + P_1$ & $C_4 + P_4 + P_2$ & $C_4 + 2P_3$ & $P_9 + P_1$
\\ \hline
$P_8 + P_2$ & $P_7 + P_3$ & $P_6 + P_4$ & $2P_5$ \\ \hline
\end{tabular}
\end{center}
Each graph $H$ in the table corresponds, up to isomorphism, to the
choice of one set of eight vertices of $G$.  Hence the question is
whether it is actually possible to extend this initial choice to the
whole graph, by coloring the remaining $17$ vertices, in such a way
that we do obtain a b-coloring with eight colors.  

We ran a computer search to examine each case separately.  It turns
out that in eleven of these cases the computer went through the
examination of the case without finding any b-coloring that extends
the position of the eight b-vertices given by graph $H$.  The
remaining case is when the pattern graph is $C_4 + P_4 + P_2$; in that
case the computer found twelve graphs.  Some of these twelve graphs
are isomorphic (by permuting rows, columns or colors), and up to
isomorphism there are only two graphs; these are shown in
Figure~\ref{fig:2K5K5}.  So these are the only two b-colorings of
$K_5\square K_5$ with $8$ colors.

Inspired by the preceding result we tried the same approach for $n=6$,
using the pattern graph $2C_4+P_2+2P_1$, and we found by hand the
graph shown on Figure~\ref{fig:K6K6}.  

\begin{figure}
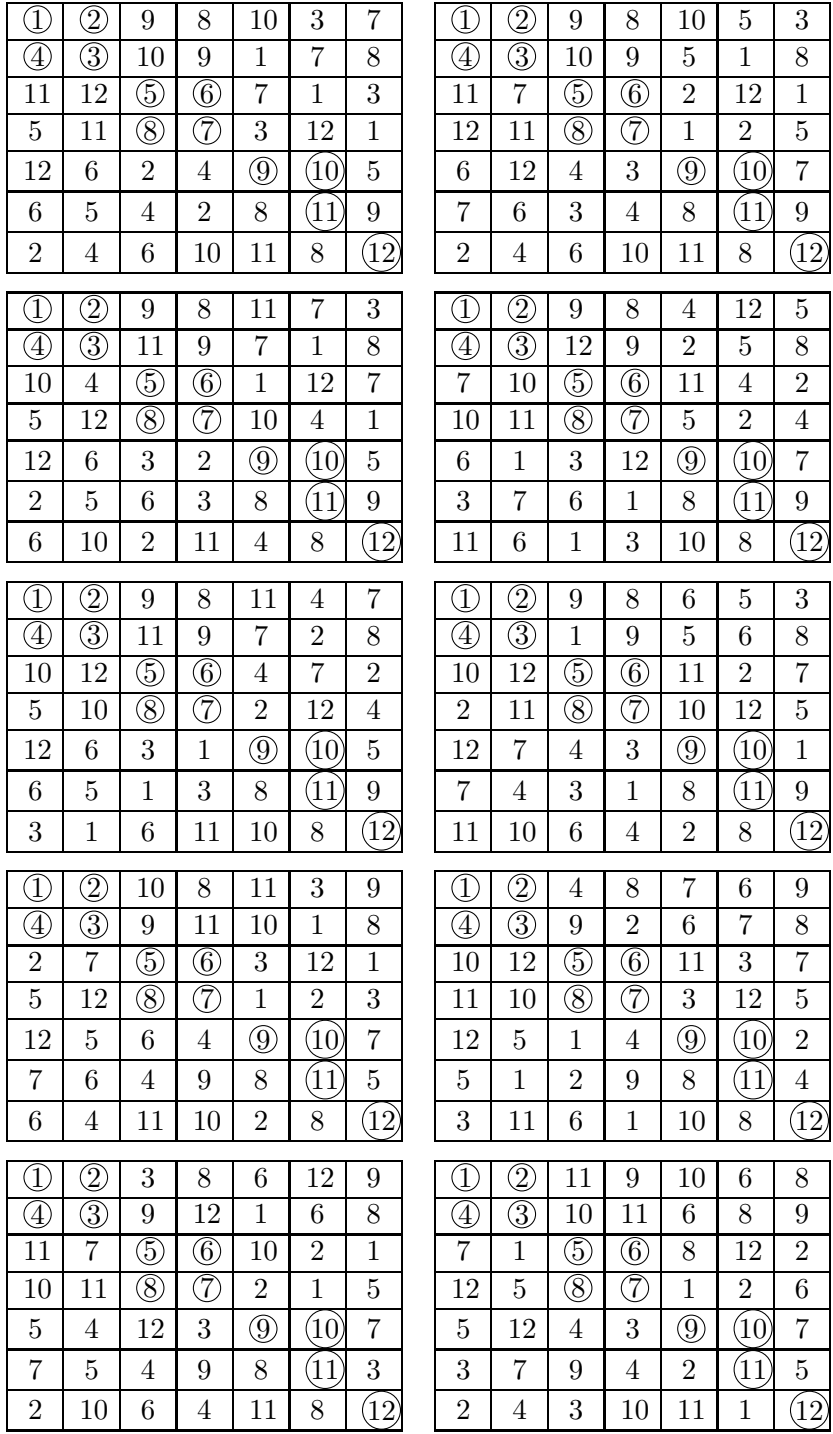

\begin{center}
\begin{tabular}{cc}
\begin{tabular}{|C|C|C|C|C|C|C|}
    \hline
\circled{1}& \circled{2}& 9& 8& 10& 3& 7 \\ \hline 
\circled{4}& \circled{3}& 10& 9& 1& 7& 8 \\ \hline  
11& 12& \circled{5}& \circled{6}& 7& 1& 3 \\ \hline  
5& 11& \circled{8}& \circled{7}& 3& 12& 1 \\ \hline  
12& 6& 2& 4& \circled{9}& \circled{10}& 5 \\ \hline  
6& 5& 4& 2& 8& \circled{11}& 9 \\ \hline  
2& 4& 6& 10& 11& 8& \circled{12} \\ \hline 
\end{tabular}
&
\begin{tabular}{|C|C|C|C|C|C|C|}
    \hline
\circled{1}& \circled{2}& 9& 8& 10& 5& 3 \\ \hline  
\circled{4}& \circled{3}& 10& 9& 5& 1& 8 \\ \hline  
11& 7& \circled{5}& \circled{6}& 2& 12& 1 \\ \hline  
12& 11& \circled{8}& \circled{7}& 1& 2& 5 \\ \hline  
6& 12& 4& 3& \circled{9}& \circled{10}& 7 \\ \hline  
7& 6& 3& 4& 8& \circled{11}& 9 \\ \hline  
2& 4& 6& 10& 11& 8& \circled{12} \\ \hline 
\end{tabular}
\end{tabular}

\medskip

\begin{tabular}{cc}
\begin{tabular}{|C|C|C|C|C|C|C|}
    \hline
\circled{1}& \circled{2}& 9& 8& 11& 7& 3 \\ \hline 
\circled{4}& \circled{3}& 11& 9& 7& 1& 8 \\ \hline 
10& 4& \circled{5}& \circled{6}& 1& 12& 7 \\ \hline 
5& 12& \circled{8}& \circled{7}& 10& 4& 1 \\ \hline 
12& 6& 3& 2& \circled{9}& \circled{10}& 5 \\ \hline 
2& 5& 6& 3& 8& \circled{11}& 9 \\ \hline 
6& 10& 2& 11& 4& 8& \circled{12} \\ \hline
\end{tabular}
&
\begin{tabular}{|C|C|C|C|C|C|C|}
    \hline
\circled{1}& \circled{2}& 9& 8& 4& 12& 5 \\ \hline 
\circled{4}& \circled{3}& 12& 9& 2& 5& 8 \\ \hline 
7& 10& \circled{5}& \circled{6}& 11& 4& 2 \\ \hline 
10& 11& \circled{8}& \circled{7}& 5& 2& 4 \\ \hline 
6& 1& 3& 12& \circled{9}& \circled{10}& 7 \\ \hline 
3& 7& 6& 1& 8& \circled{11}& 9 \\ \hline 
11& 6& 1& 3& 10& 8& \circled{12} \\ \hline
\end{tabular}
\end{tabular}

\medskip

\begin{tabular}{cc}
\begin{tabular}{|C|C|C|C|C|C|C|}
    \hline
\circled{1}& \circled{2}& 9& 8& 11& 4& 7 \\ \hline 
\circled{4}& \circled{3}& 11& 9& 7& 2& 8 \\ \hline 
10& 12& \circled{5}& \circled{6}& 4& 7& 2 \\ \hline 
5& 10& \circled{8}& \circled{7}& 2& 12& 4 \\ \hline
12& 6& 3& 1& \circled{9}& \circled{10}& 5 \\ \hline 
6& 5& 1& 3& 8& \circled{11}& 9 \\ \hline 
3& 1& 6& 11& 10& 8& \circled{12} \\ \hline
\end{tabular}
&
\begin{tabular}{|C|C|C|C|C|C|C|}
    \hline
\circled{1}& \circled{2}& {9\ }& 8& 6& 5& 3 \\ \hline 
\circled{4}& \circled{3}& 1& 9& 5& 6& 8 \\ \hline 
10& 12& \circled{5}& \circled{6}& 11& 2& 7 \\ \hline 
2& 11& \circled{8}& \circled{7}& 10& 12& 5 \\ \hline 
12& 7& 4& 3& \circled{9}& \circled{10}& 1 \\ \hline 
7& 4& 3& 1& 8& \circled{11}& 9 \\ \hline 
11& 10& 6& 4& 2& 8& \circled{12} \\ \hline
\end{tabular}
\end{tabular}

\medskip

\begin{tabular}{cc}
\begin{tabular}{|C|C|C|C|C|C|C|}
    \hline
\circled{1}& \circled{2}& 10& 8& 11& 3& 9 \\ \hline 
\circled{4}& \circled{3}& 9& 11& 10& 1& 8 \\ \hline 
2& 7& \circled{5}& \circled{6}& 3& 12& 1 \\ \hline 
5& 12& \circled{8}& \circled{7}& 1& 2& 3 \\ \hline 
12& 5& 6& 4& \circled{9}& \circled{10}& 7 \\ \hline
7& 6& 4& 9& 8& \circled{11}& 5 \\ \hline 
6& 4& 11& 10& 2& 8& \circled{12} \\ \hline
\end{tabular}
&
\begin{tabular}{|C|C|C|C|C|C|C|}
    \hline
\circled{1}& \circled{2}& 4& 8& 7& 6& 9 \\ \hline
\circled{4}& \circled{3}& 9& 2& 6& 7& 8 \\ \hline 
10& 12& \circled{5}& \circled{6}& 11& 3& 7 \\ \hline 
11& 10& \circled{8}& \circled{7}& 3& 12& 5 \\ \hline 
12& 5& 1& 4& \circled{9}& \circled{10}& 2 \\ \hline 
5& 1& 2& 9& 8& \circled{11}& 4 \\ \hline 
3& 11& 6& 1& 10& 8& \circled{12} \\ \hline
\end{tabular}
\end{tabular}

\medskip

\begin{tabular}{cc}
\begin{tabular}{|C|C|C|C|C|C|C|}
    \hline
\circled{1}& \circled{2}& 3& 8& 6& 12& 9 \\ \hline 
\circled{4}& \circled{3}& 9& 12& 1& 6& 8 \\ \hline 
11& 7& \circled{5}& \circled{6}& 10& 2& 1 \\ \hline 
10& 11& \circled{8}& \circled{7}& 2& 1& 5 \\ \hline 
5& 4& 12& 3& \circled{9}& \circled{10}& 7 \\ \hline 
7& 5& 4& 9& 8& \circled{11}& 3 \\ \hline 
2& 10& 6& 4& 11& 8& \circled{12} \\ \hline 
\end{tabular}
&
\begin{tabular}{|C|C|C|C|C|C|C|}
    \hline
\circled{1}& \circled{2}& 11& 9& 10& 6& 8 \\ \hline 
\circled{4}& \circled{3}& 10& 11& 6& 8& 9 \\ \hline
7& 1& \circled{5}& \circled{6}& 8& 12& 2 \\ \hline 
12& 5& \circled{8}& \circled{7}& 1& 2& 6 \\ \hline 
5& 12& 4& 3& \circled{9}& \circled{10}& 7 \\ \hline 
3& 7& 9& 4& 2& \circled{11}& 5 \\ \hline 
2& 4& 3& 10& 11& 1& \circled{12} \\ \hline
\end{tabular}
\end{tabular}
\end{center}
\caption{Some b-colorings of $K_7\square K_7$ with $12$ colors}
\label{fig:K7K7}
\end{figure}

For $n=7$ we started from the pattern graph $2C_4+P_4+P_2$.  This was
a fruitful choice and the computer produced tens of thousands of
solutions, namely b-colorings of $K_7\square K_7$ with $12$ colors (it
was not possible to test isomorphism between them); we show a few of
them of them in Figure~\ref{fig:K7K7}.  It might be that for $n=6$ or
$7$ there are solutions with a different pattern graph.  (Moreover,
when $n\ge 6$, the pattern graph does not necessarily have maximum
degree~$2$.)  

We notice that in all the solutions we found each color class has only
one b-vertex.

For $n=8$ or $9$ the number of cases to explore became too high, even
when we started, for $n=9$, from the pattern graph $3C_4+P_4+P_2$
which should be a good candidate.

\begin{figure}
\begin{center}
\begin{tabular}{cc}
\begin{tabular}{|c|c|c|c|c|}
    \hline
\circled{1}&	\circled{2}&	8&	4&	6 \\ \hline	
\circled{4}&	\circled{3}&	6&	8&	2 \\ \hline	
3&	7&	\circled{5}&	\circled{6}&	1 \\ \hline	
5&	1&	2&	\circled{7}&	3 \\ \hline	
7&	5&	4&	2&	\circled{8} \\ \hline
\end{tabular}
&
\begin{tabular}{|c|c|c|c|c|}
\hline
\circled{1}&	\circled{2}&	8&	4&	6 \\ \hline	
\circled{4}&	\circled{3}&	6&	8&	2 \\ \hline	
7&	1&	\circled{5}&	\circled{6}&	3 \\ \hline	
3&	5&	2&	\circled{7}&	1 \\ \hline	
5&	7&	4&	2&	\circled{8} \\ \hline	
\end{tabular}
\end{tabular}
\end{center} 
\caption{b-colorings of $K_5\square K_5$ with $8$
colors. The b-vertices are circled.}\label{fig:2K5K5}
\end{figure}

\bigskip

\begin{figure}
\begin{center}
\begin{tabular}{|C|C|C|C|C|C|}
    \hline
\circled{1}& \circled{2}& 4& 9& 5& 7 \\ \hline 
\circled{4}& \circled{3}& 9& 2& 7& 5 \\ \hline  
10& 8& \circled{5}& \circled{6}& 3& 1 \\ \hline  
6& 10& \circled{8}& \circled{7}& 1& 3 \\ \hline  
8& 6& 2& 4& \circled{9}& \circled{10} \\ \hline  
3& 1& 7& 5& 10& 9 \\ \hline  
\end{tabular}
\end{center}
\caption{A b-coloring of $K_6\square K_6$ with $10$ colors}
\label{fig:K6K6}
\end{figure}


\clearpage
\section*{Appendix: Proof of Lemma~\ref{lem:beta}}

For $c=1, \ldots, 8$ let $S_c$ be the set of vertices of color~$c$ in
$f$.  For any set $X\subseteq V(G)$, we write $f(X)=\{f(x) \mid x\in
X\}$.

We will frequently use the following argument.  Every vertex of $G$
has degree~$8$.  If $x$ is a b-vertex, it has a neighbor of each of
the seven colors different from $f(x)$, and the remaining neighbor of
$x$ has one of these colors, so $x$ has exactly one pair of neighbors
of the same color (a \emph{repeat}).  So if any vertex has two repeats
in its neighborhood it cannot be a b-vertex.

Suppose that the lemma does not hold.  So, up to symmetry, we have
$\beta_1\ge 3$.  So we may assume, up to symmetry, that $v_{1,c}$ is a
b-vertex of color~$c$ for $c=1,2,3$ and $v_{1,j}\in S_j$ for $j=4,5$.
It must be that colors $6,7$ and $8$ appear in column~$c$ for each
$c=1,2,3$, and consequently there is a unique integer $t_c$ from
$\{2,3,4,5\}$ such that the color of $v_{t_c,c}$ is not from
$\{c,6,7,8\}$.  Up to symmetry, this leads to the following three
cases (See Figure~\ref{fi$G_1$23}).

\begin{figure}[h]
\hspace{.3em}
\begin{tabular}{|c|c|c|c|c|}
    \hline
\circled{1} & \circled{2} & \circled{3} & 4 & 5 \\
\hline
6  & 7  &  8 &   &   \\
\hline
7  & 8  & 6  &   &   \\
\hline
8  & 6  &  7 &   &   \\
\hline
  &   &   &   &   \\
\hline
\end{tabular}
\hfill
\begin{tabular}{|c|c|c|c|c|}
    \hline
\circled{1} & \circled{2} & \circled{3} & 4 & 5 \\
\hline
6  & 7  &    &   &   \\
\hline
7  & 8  & 6  &   &   \\
\hline
8  & 6  &  7 &   &   \\
\hline
  &   & 8  &   &   \\
\hline
\end{tabular}
\hfill
\begin{tabular}{|c|c|c|c|c|}
    \hline
\circled{1} & \circled{2} & \circled{3} & 4 & 5 \\
\hline
6  & 7  &  8  &   &   \\
\hline
7  & 6  &    &   &   \\
\hline
8  &    &  6 &   &   \\
\hline
  & 8  & 7  &   &   \\
\hline
\end{tabular}
\hspace{.3em}
\caption{Cases 1, 2 and 3}\label{fi$G_1$23}
\end{figure}

\medskip

\noindent{\it Case 1: $t_1=t_2=t_3$.} We may assume that $t_1=5$ and
that $v_{2,1}, v_{3,3}, v_{4,2} \in S_6$, $v_{2,2}, v_{3,1}, v_{4,3}
\in S_7$, and $v_{2,3}, v_{3,2}, v_{4,1}\in S_8$.  See
Figure~\ref{fi$G_1$23}, left.  Since $v_{2,1}$ has two repeats in its
neighborhood (with colors $7$ and $8$), it is not a b-vertex.
Likewise, $v_{i,j}$ is not a b-vertex for each $i=2,3,4$ and $j=
1,2,3$.  Apart from these nine vertices, only vertices $v_{5,4}$ and
$v_{5,5}$ may have colors $6,7,8$, so $\{b_6,b_7,b_8\} \subseteq
\{v_{5,4}, v_{5,5}\}$, which is impossible.

\medskip

\noindent{\it Case 2: $t_1=t_2$ and $t_1\neq t_3$.} We may assume that
$t_1=t_2=5$ and $t_3=2$, and that $v_{2,1}, v_{3,3}, v_{4,2} \in S_6$,
$v_{2,2}, v_{3,1}, v_{4,3} \in S_7$, and $v_{3,2}, v_{4,1}, v_{5,3}
\in S_8$.  See Figure~\ref{fi$G_1$23}, center.  For $i=3,4$ and
$j=1,2,3$, vertex $v_{i,j}$ has two repeats in its neighborhood, so it
is not a b-vertex.  Likewise $v_{2,3}$ is not a b-vertex.  Therefore:
\begin{equation}\label{eq:b6b7}
b_6\in \{v_{2,1}, v_{5,4}, v_{5,5}\} \mbox{ and } b_7\in \{v_{2,2},
v_{5,4}, v_{5,5}\}.
\end{equation}
The only neighbor of $v_{1,4}$ that can have color $6$ or $7$ is
$v_{5,4}$, and it cannot have both colors, so $b_4\neq v_{1,4}$.  By
the same argument, $b_4\neq v_{5,5}$.  So $b_4$ is either in
$\{v_{5,1}, v_{5,2}\}$ or in $\{v_{2,5}, v_{3,5}, v_{4,5}\}$.

Suppose that $b_4=v_{5,1}$.  Then $f(\{v_{5,2},v_{5,4},v_{5,5}\})=$
$\{2,3,5\}$, and it follows from (\ref{eq:b6b7}) that $b_6=v_{2,1}$
and $b_7=v_{2,2}$.  But then colors $1,2,3,4,5$ should appear in
$\{v_{2,3}, v_{2,4}, v_{2,5}\}$, which is impossible.  Hence $b_4\neq
v_{5,1}$, and similarly $b_4\neq v_{5,2}$.  So $b_4$ is in $\{v_{2,5},
v_{3,5}, v_{4,5}\}$.  Likewise $b_5$ is in $\{v_{2,4}, v_{3,4},
v_{4,4}\}$.

Suppose that $b_4=v_{2,5}$.  Then $v_{2,4}\in S_8$, for otherwise
$b_4$ has no neighbor of color $8$, and so $b_5\in \{v_{3,4},
v_{4,4}\}$.  Now $v_{2,1}$ is not a b-vertex, because it has two
repeats (with colors $7$ and $8$), and similarly $v_{2,2}$ is not a
b-vertex.  So, by (\ref{eq:b6b7}), we have $b_6, b_7\in \{v_{5,4},
v_{5,5}\}$.  But then $b_5$ has two repeats (with color $8$ and either
$6$ or $7$), a contradiction.  So $b_4\in \{v_{3,5}, v_{4,5}\}$.
Likewise, $b_5\in\{v_{3,4}, v_{4,4}\}$.  Hence colors $4$ and $5$ do
not appear in $\{v_{2,4}, v_{2,5}\}$.  Whichever way colors $4$ and
$5$ appear in $\{v_{2,3}, v_{5,1}, v_{5,2}\}$, one of $b_6,b_7$ has no
neighbor of some color $j\in\{4,5\}$, so it is not a b-vertex.  So, by
(\ref{eq:b6b7}) and up to symmetry, we may assume that $b_6=v_{5,4}$.
Hence $v_{5,5}$ has color $7$, for otherwise $b_6$ has no neighbor of
color $7$.  Now we have $f(\{v_{2,5}, v_{3,4}, v_{4,5}\})=$
$\{1,2,3\}$ (for $b_4$ to be a b-vertex) and similarly $f(\{v_{2,4},
v_{3,4}, v_{4,5}\})=$ $\{1,2,3\}$ (for $b_5$); but this entails that
$v_{2,4}$ and $v_{2,5}$ have the same color, a contradiction.

\medskip

\noindent{\it Case 3: $t_1, t_2$ and $t_3$ are pairwise different.} We
may assume that $t_1=5$, $t_2=4$ and $t_3=3$, and that $v_{2,1},
v_{3,2}, v_{4,3}\in S_6$, $v_{2,2}, v_{3,1}, v_{5,3}\in S_7$, and
$v_{2,3}, v_{4,1}, v_{5,2}\in S_8$.  See Figure~\ref{fi$G_1$23},
right.  Vertices $v_{2,1}$, $v_{2,2}$ and $v_{2,3}$ are not b-vertices
because they have two repeats.  Likewise, $v_{5,1}, v_{4,2}$ and
$v_{3,3}$ are not b-vertices.  Hence $b_4$ is in column $5$ and $b_5$
is in column $4$.

Suppose that $b_4=v_{2,5}$ and $b_5=v_{2,4}$.  Then $f(\{v_{3,4},
v_{4,4}, v_{5,4}\})=$ $\{1,2,3\}$ $f(\{v_{3,5}, v_{4,5}, v_{5,5}\})=$
$\{1,2,3\}$.  Rows $3,4,5$ play symmetric roles, so we may assume that
$b_7=v_{3,1}$.  So $f(\{v_{3,3}, v_{5,1}\})=$ $\{4,5\}$ and
$f(\{v_{3,4}, v_{3,5}\})=$ $\{2,3\}$.  Then $v_{3,2}$ is not a
b-vertex because it has two repeats (with colors $2$ and $7$), so
$b_6=v_{4,3}$.  Similarly, $b_8=v_{5,2}$.  Now, whichever way colors
$4$ and $5$ appear in $\{v_{5,1}, v_{4,2}, v_{3,3}\}$, one of $b_6,
b_7, b_8$ has no neighbor of some color $j\in\{4,5\}$, a
contradiction.  

Therefore we may assume, up to symmetry, that $b_4=v_{3,5}$.  Then
$v_{3,4}$ has color $8$, for otherwise $b_4$ has no neighbor of color
$8$.  Then $v_{3,1}$ and $v_{3,2}$ are not b-vertices because they
have two repeats.  If $b_5=v_{4,4}$, then by a similar argument
$v_{4,5}$ has color $7$ and $v_{4,1}, v_{4,3}$ are not b-vertices, and
only colors $1,2,3$ appear in $\{v_{2,4}, v_{2,5}, v_{3,3}, v_{4,2},
v_{5,4}, v_{5,5}\}$, so no vertex can play the role of $b_6$, a
contradiction.  Hence $b_5\neq v_{4,4}$, and similarly $b_5\neq
v_{4,5}$, so $b_5= v_{4,2}$.  We must have $f(\{v_{2,5}, v_{4,4},
v_{5,4}\})=\{1,2,3\}$ for $b_5$ to be a b-vertex.  The same argument
for $b_4$ implies that some vertex $v\in\{v_{4,5}, v_{5,5}\}$
satisfies $f(v)\in\{1,2,3\}$, say (by symmetry of rows $4$ and $5$)
$v=v_{4,5}$.  Then $S_7=\{v_{2,2}, v_{3,1}, v_{5,3}\}$, so $b_7=
v_{5,3}$.  We have $f(\{v_{5,4}, v_{5,5}\})=$ $\{1,2\}$, $f(v_{3,3})=
5$ and $f(v_{5,1})=4$.  Now the symmetry of rows $4,5$ is restored, so
we also have $b_6=v_{4,3}$ and $f(\{v_{4,4}, v_{4,5}\})=$ $\{1,2\}$
and $f(v_{4,2})=4$.  In $S_8$ every vertex has two repeats, so there
is no b-vertex in $S_8$, a contradiction.  This completes the proof of
Lemma~\ref{lem:beta}.  \hfill $\Box$

\medskip

The next lemma could also be used to help decrease the number of
subcases to explore.

\begin{lemma}\label{lem:size}
    The sets $S_1, \ldots, S_8$ satisfy the following properties: 
\vspace{-.2cm}
\begin{itemize}
\setlength{\itemsep}{-.1cm}
\item
One of $S_1, \ldots, S_8$ has size~$4$, and the other seven sets
have size~$3$.  
\item
If $|S_i|=3$, then $S_i\setminus \{b_i\}$ contains a vertex that
has at least three neighbors in $B\setminus(\{b_i\}\cup N(b_i))$.
\end{itemize}
\end{lemma}
\noindent{\it Proof.} For each $i\in\{1, \ldots, 8\}$, by
Lemma~\ref{lem:beta} the row that contains $b_i$ contains at most one
element of $B \setminus \{b_i\}$, and the same holds for the column
that contains $b_i$.  Hence at least five members of
$B\setminus\{b_i\}$ have a neighbor in $S_i\setminus\{b_i\}$.  So if
$S_i\setminus\{b_i\}$ contains only one vertex $x_i$, then three of
these five elements are either on the row of $x_i$ or on the column of
$x_i$, which violates Lemma~\ref{lem:beta}.  So $|S_i|\ge 3$, and if
$|S_i|= 3$ then the second item of the lemma must hold.  Also we have
$25 = \sum_{i=1}^8 |S_i| \ge 8 \times 3 = 24$, which implies the first
item of the lemma.  \hfill $\Box$



\begin{thebibliography}{9}
    
\bibitem{F}
T.~Faik.  {\it La b-continuit\'e des b-colorations: Complexit\'e,
propri\'et\'es structurelles et algorithmes.} Ph.D. Thesis, Univ.
Orsay, France, 2005.

\bibitem{IM}
R.W.~Irving, D.F.~Manlove.  The b-chromatic number of a graph.  {\it
Discrete Applied Mathematics} 91 (1999) 127--141.

\bibitem{JavOm}
R.~Javadi, B.~Omoomi.  On b-coloring of cartesian product of graphs.
{\it Ars Combinatoria} 107 (2012) 521--536.
\end{thebibliography}
\end{document}